\documentclass[11pt]{amsart}
\usepackage{amssymb,hyperref}

\newtheorem{theorem}{Theorem}

\newtheorem{conjecture}[theorem]{Conjecture}
\newtheorem{corollary}[theorem]{Corollary}
\newtheorem{definition}[theorem]{Definition}
\newtheorem{example}[theorem]{Example}

\newtheorem{fact}[theorem]{Fact}
\newtheorem{lemma}[theorem]{Lemma}
\newtheorem{problem}[theorem]{Problem}
\newtheorem{proposition}[theorem]{Proposition}
\newtheorem{question}[theorem]{Question}
\newtheorem{remark}[theorem]{Remark}

\newtheorem*{theorem*}{Theorem}
\newtheorem*{proposition*}{Proposition}

\newcommand{\bcon}{\begin{conjecture}}
\newcommand{\econ}{\end{conjecture}}
\newcommand{\bcor}{\begin{corollary}}
\newcommand{\ecor}{\end{corollary}}
\newcommand{\bdf}{\begin{definition}}
\newcommand{\edf}{\end{definition}}
\newcommand{\beq}{\begin{equation}}
\newcommand{\eeq}{\end{equation}}
\newcommand{\bexa}{\begin{example}}
\newcommand{\eexa}{\end{example}}
\newcommand{\bfac}{\begin{fact}}
\newcommand{\efac}{\end{fact}}
\newcommand{\blem}{\begin{lemma}}
\newcommand{\elem}{\end{lemma}}
\newcommand{\bprb}{\begin{problem}}
\newcommand{\eprb}{\end{problem}}
\newcommand{\bpro}{\begin{proposition}}
\newcommand{\epro}{\end{proposition}}

\newcommand{\bque}{\begin{question}}
\newcommand{\eque}{\end{question}}
\newcommand{\brem}{\begin{remark}}
\newcommand{\erem}{\end{remark}}
\newcommand{\bthm}{\begin{theorem}}
\newcommand{\ethm}{\end{theorem}}
\newcommand{\bmat}{\begin{matrix}}
\newcommand{\emat}{\end{matrix}}

\newcommand{\be}{\begin}
\newcommand{\en}{\end}
\newcommand{\bpr}{\begin{proof}}
\newcommand{\epr}{\end{proof}}

\newcommand{\lb}{\label}
\newcommand{\la}{\langle}
\newcommand{\ra}{\rangle}
\newcommand{\comment}[1]{\,}

\newcommand{\cal}{\mathcal}

\newcommand{\Z}{\mathbb Z}

\newcommand{\C}{\mathbb C}

\newcommand{\T}{\mathbb T}

\newcommand{\ve}{\varepsilon}

\newcommand{\FT}{{\cal F\cal T}}
\renewcommand{\T}{{\cal T}}

\title{$G$-Character varieties for $G=\operatorname{SO}(n,\C)$ and other not simply connected groups}

\author{Adam S. Sikora}
\address{244 Math Bldg, University at Buffalo, SUNY, Buffalo, NY 14260}
\email{asikora@buffalo.edu}
\thanks{The author acknowledges support from U.S. National Science Foundation grants DMS 1107452, 1107263, 1107367 "RNMS: GEometric structures And Representation varieties" (the GEAR Network).}


\begin{document}

\thispagestyle{empty}

\begin{abstract} We describe the relation between $G$-character varieties, $X_G(\Gamma),$ and $G/H$-character varieties, where $H$ is a finite, central subgroup of $G$. In particular, we find finite generating sets for $\C[X_{G/H}(\Gamma)]$ for classical groups $G$ and $H$ as above. By applying this approach to $SO(4,\C)=(SL(2,\C)\times SL(2,\C))/\Z/2$ we find an explicit description of $\C[X_{SO(4,\C)}(F_2)]$ for the free group on two generators, $F_2.$

In the second part of the paper, we prove several properties of
$SO(2n,\C)$-character varieties. This is a particularly interesting class of
character varieties because unlike for all other classical groups $G$, the
coordinate rings $\C[X_{G}(\Gamma)]$ are generally not generated by the trace
functions, $\tau_\gamma$, for $\gamma\in \Gamma,$ for $G=SO(2n,\C)$. In fact, we prove that the coordinate ring $\C[X_{SO(2n,\C)}(\Gamma)]$ is not even generated by ``generalized trace functions,'' $\tau_{\gamma,V},$ for all $\gamma\in \Gamma$ and all representations $V$ of $SO(2n,\C)$ for $n=2$ and groups $\Gamma$ of corank $\geq 2$.
\end{abstract}

\maketitle

\pagestyle{myheadings}
\markboth{\hfil{\sc ADAM S. SIKORA}\hfil}
{\hfil{\sc $G$-CHARACTER VARIETIES FOR $G=\operatorname{SO}(n,\C)$}\hfil}

%
\section{Introduction}
%

Let $G$ be an affine reductive algebraic group over $\C.$\footnote{Throughout the paper, the field of complex numbers can be replaced an arbitrary algebraically closed field of characteristic zero.} For every group $\Gamma$ generated by some $\gamma_1,...,\gamma_N$, the space of all $G$-representations of $\Gamma$ forms an algebraic subset, $Hom(\Gamma,G),$ of $G^N$ on which $G$ acts by conjugating representations.
The categorical quotient of that action $$X_G(\Gamma)=Hom(\Gamma,G)//G$$ is the $G$-character variety of $\Gamma,$ cf. \cite{LM, S-char} and the references within.

In this paper, we study character varieties for groups which are not simply connected.
In Sections \ref{s_G-G'}-\ref{s_gen} we describe the relations between $X_{G}(\Gamma)$ and $X_{G/H}(\Gamma),$ for all finite, central subgroups $H$ of $G$. In particular, we find finite generating sets for $\C[X_{G/H}(\Gamma)]$ for classical groups $G$ and $H$ as above, cf. Theorem \ref{generators-G'}, Corollary \ref{gen-M'}, and Proposition \ref{m^N}. By applying this approach to $SO(4,\C)=(SL(2,\C)\times SL(2,\C))/\Z/2$ we find an explicit description of $\C[X_{SO(4,\C)}(F_2)]$ for the free group on two generators, $F_2,$ cf. Proposition \ref{SO4-desc}.

In Sections \ref{s_trace_alg}-\ref{s_Coord-gen2}, we use the above results to prove several statements about the $SO(2n,\C)$-character varieties. This is a particularly interesting class of character varieties because unlike for all other classical groups $G$, the coordinate ring $\C[X_{G}(\Gamma)]$ does not always coincide with the $G$-trace algebra of $\Gamma,$ $\T_G(\Gamma),$ for $G=SO(2n,\C),$ cf. \cite{S-gen}.

Let us recall the definition of trace algebras. For a representation $V$ of $G$ and $\gamma\in \Gamma$
let $$\tau_{\gamma,V}([\rho])=tr(\phi_V\rho(\gamma)),$$
where $\phi_V$ is the homomorphism $\Gamma\to GL(V)$ induced by $V.$
Then $\tau_{\gamma,V}$ is a regular function on $X_{G}(\Gamma).$ For a matrix group $G\subset GL(n,\C)$, the $G$-trace algebra, $\T_G(\Gamma)$,
is the $\C$-subalgebra of $\C[X_{G}(\Gamma)]$ generated by $\tau_\gamma=\tau_{\gamma,\C^n},$  for all $\gamma\in \Gamma.$

Let us also recall from \cite{S-gen} that the full $G$-trace algebra, $\FT_G(\Gamma)$,
is the $\C$-subalgebra of $\C[X_{G}(\Gamma)]$ generated by $\tau_{\gamma,V},$ for all representations $V$ of $G$ and for all $\gamma\in \Gamma.$ Addressing a question posted in \cite{S-gen}, we investigate the extensions $$\T_{SO(2n,\C)}(\Gamma) \subset \FT_{SO(2n,\C)}(\Gamma) \subset \C[X_{SO(2n,\C)}(\Gamma)].$$

By \cite[Thm 8]{S-gen}\footnote{A partial version of this result (with infinite generating sets) can be also found in \cite{ATZ}}, $\C[X_{SO(2n,\C)}(\Gamma)]$ is a $\T_{SO(2n,\C)}(\Gamma)$-algebra generated by $Q_{2n}(g_1,...,g_n)$ for all words $g_1,...,g_n$ in $\gamma_1,...,\gamma_N$ of length at most $\nu_n-1$ in which the number of inverses is not larger than half of the length of the word. (See  \cite{S-gen} for the definition of $\nu_n$ and of $Q_{2n}(g_1,...,g_n)$. We recall the definition of the latter in Sec. \ref{s_gen}.) Here is a stronger version of that result for free groups:

\begin{proposition}(Proof in Sec. \ref{s_trace_alg}.)\lb{module_gen}
If $\Gamma$ is free on $\gamma_1,...,\gamma_N,$ then $\C[X_{SO(2n,\C)}(\Gamma)]$ is a $\T_{SO(2n,\C)}(\Gamma)$-module generated by $1$ and by $Q_{2n}(g_1,...,g_n)$ for words $g_1,...,g_n\in \Gamma$ in $\gamma_1,...,\gamma_N$ of length at most $\nu_n-1$ in which the number of inverses is not larger than the half of the length of the word.
\end{proposition}

The determination of whether the extension
$$\FT_{SO(2n,\C)}(\Gamma)\subset \C[X_{SO(2n,\C)}(\Gamma)]$$
is proper is a delicate problem.
We analyze this extension in detail for $\Gamma=F_2$ and $n=2$ in Sections \ref{s_trace_alg}-\ref{s_Coord-gen2}. In particular we prove:

\begin{theorem}\lb{Coord-gen} (Proof in Sec. \ref{s_Coord-gen1} and \ref{s_Coord-gen2}.)
\begin{enumerate}
\item $\C[X_{SO(4,\C)}(F_2)]$ is a $\T_{SO(4,\C)}(F_2)$-module generated by
$$1,Q_4(\gamma_i,\gamma_j),\ \text{for}\ 1\leq i\leq j\leq 2,$$
and by
$$Q_4(\gamma_1\gamma_2,\gamma_1\gamma_2),
Q_4(\gamma_1\gamma_2^{-1},\gamma_1\gamma_2^{-1}),
Q_4(\gamma_1\gamma_2^{-1},\gamma_2),
Q_4(\gamma_2\gamma_1^{-1},\gamma_1).$$
\item $\C[X_{SO(4,\C)}(F_2)]$ is a $\FT_{SO(4,\C)} (F_2)$-module generated by
$$1, Q_4(\gamma_1,\gamma_2), Q_4(\gamma_1\gamma_2^{-1},\gamma_2),\ \text{and}\ Q_4(\gamma_2\gamma_1^{-1},\gamma_1).$$
\item None of the generators of part (2) can be expressed as a linear combination of others with coefficients in $\FT_{SO(4,\C)}(F_2)$. Consequently,\\
$\FT_{SO(4,\C)} (F_2)\varsubsetneq \C[X_{SO(4,\C)}(F_2)].$
\end{enumerate}
\end{theorem}

Since every epimorphism $\phi:\Gamma\to F_2$ induces an epimorphism
$$\phi_*:\C[X_{SO(4,\C)}(\Gamma)]\to \C[X_{SO(4,\C)}(F_2)]$$ 
mapping $\FT_{SO(4,\C)}(\Gamma)$ to $\FT_{SO(4,\C)}(F_2),$ we have

\bcor 
For every group $\Gamma$ of corank $\geq 2$ (i.e. having $F_2$ as a quotient),
$\FT_{SO(4,\C)}(\Gamma)$ is a proper subalgebra of $\C[X_{SO(4,\C)}(\Gamma)]$.
\ecor

\bcon 
For every $n>2$ and every group $\Gamma$ of corank $\geq 2,$
$\FT_{SO(2n,\C)}(\Gamma)$ is a proper subalgebra of $\C[X_{SO(2n,\C)}(\Gamma)]$.
\econ

\noindent {\bf Acknowledgments}
We would like to thank the anonymous referee for careful reading of this paper and helpful comments.

%
\section{The relation between $X_G(\Gamma)$ and $X_{G/H}(\Gamma)$}
\lb{s_G-G'}
%

Let $Hom^0(\Gamma,G)$ denote the connected component of the trivial representation in $Hom(\Gamma,G).$
Let $X_G^0(\Gamma)$ denote its image in $X_G(\Gamma)$.

Let $H$ be a finite, central subgroup of $G.$
Then $Hom(\Gamma,H)$ acts on $Hom(\Gamma,G)$ by multiplication,
$$\sigma\cdot \rho(\gamma)=\sigma(\gamma)\cdot \rho(\gamma),$$
for $\gamma\in \Gamma.$
Let $Hom'(\Gamma,H)$ be the set of elements of $Hom(\Gamma,H)$ which map $Hom^0(\Gamma,G)$ to itself.

\begin{proposition}\label{G-G'} (1) The projection $G\to G/H$ induces an isomorphism
$$\phi: Hom^0(\Gamma,G)/Hom'(\Gamma,H)\to Hom^0(\Gamma,G/H).$$
(2) Furthermore, $\phi$ descends to an isomorphism
$$\psi: X^0_G(\Gamma)/Hom'(\Gamma,H)\to X^0_{G/H}(\Gamma).$$
\end{proposition}

Note that since $Hom'(\Gamma,H)$ is finite, the above quotients are both ``set theory" and categorical quotients.

\noindent{\it Proof of Proposition \ref{G-G'}:}
(1) Since the extension
$$\{ e\}\to H\to G\to G/H\to \{ e\}$$
is central, it defines an element $\alpha\in H^2(G/H,H)$ such that
$f:\Gamma\to G/H$ lifts to $\tilde f: \Gamma\to G$ if and only if $f^*(\alpha)=0$ in $H^2(\Gamma,H),$ cf.
\cite[Sec. 2]{GM}.
Since $H^2(\Gamma,H)$ is discrete, the property of $f$ being ``liftable" is locally constant on $Hom(\Gamma,G)$ (in complex topology). Since the trivial representation is liftable, the above argument shows that $G\twoheadrightarrow G/H$ induces an epimorphism $Hom^0(\Gamma,G)\to Hom^0(\Gamma,G/H)$. That epimorphism factors to
$$\phi: Hom^0(\Gamma,G)/Hom'(\Gamma,H)\to Hom^0(\Gamma,G/H).$$
Since any two representations in $Hom^0(\Gamma,G)$ projecting to the same representation in $Hom^0(\Gamma,G/H)$ differ by an element of $Hom'(\Gamma,H)$, $\phi$ is $1$-$1$.

(2) follows from the fact that the action of $Hom'(\Gamma,H)$ on $Hom^0(\Gamma,G)$ commutes with the $G$-action by conjugation.
\qed

\bpro\lb{G-G'-free} The projection $G\to G/H$ induces an isomorphism\\
$\psi: X_G(F_N)/H^N\to X_{G/H}(F_N)$ for free groups, $F_N.$
\epro

\bpr For $G$ connected, $Hom(F_N,G)=G^N$ is connected and, hence, the statement is an immediate consequence of
Proposition \ref{G-G'}. Observe, however, that the only reason for considering a specific connected component of $Hom(\Gamma,G/H)$ was to make sure that representations $\Gamma\to G/H$ are liftable to $G$.
Since all representations are liftable for $\Gamma$ free, the proof above implies the statement of this proposition as well.
\epr

%
\section{Finite generating sets for $\C[X_{G/H}(\Gamma)]$}
\lb{s_gen}
%

Let us assume that $H$ is a cyclic group of order $m$ of scalar matrices in a matrix group $G.$
We will apply the above results to provide finite generating sets for the coordinate rings of $G/H$-character varieties. Since for every $\Gamma$ the natural projection
$$\pi: F_N=\la \gamma_1,...,\gamma_N\ra \to \Gamma$$
induces an epimorphism
$$\pi_*:\C[X_{G/H}(F_N)]\to \C[X_{G/H}(\Gamma)],$$
we will construct finite generating sets for $\C[X_{G/H}(\Gamma)]$ for free groups $\Gamma$ only.

Given $\gamma\in F_N=\la \gamma_1,...,\gamma_N\ra,$ let $v(\gamma)$ be a vector in $(\Z/m)^N,$ whose $i$-th component is the sum (mod $m$) of exponents of $\gamma_i$ in $\gamma$.

\begin{theorem} \label{generators-G'}
\textup{(1)} If $g_1,...,g_k$ are elements of $F_N$ such that
$\sum_{i=1}^k v(g_i)=0$, then there exists a unique $\lambda_{g_1,...,g_k}: X_{G/H}(F_N)\to \C$ such that\\ $\tau_{g_1}\cdot ...\cdot \tau_{g_k}=\lambda_{g_1,...,g_k}\psi,$ where $\psi:X_{G}(F_N)\to X_{G/H}(F_N)$ is the map of Proposition \ref{G-G'-free}.\\
\textup{(2)} Suppose that $\C[X_G(F_N)]$ is generated by $\tau_g,$ for $g$'s in a set $B\subset F_N.$ Let $\cal M$ be the set of all $\lambda_{g_1,...,g_k}$ for $g_1,...,g_k\in B$ such that $\sum_{i=1}^k v(g_k)=0$.
Then $\C[X_{G/H}(F_N)]$ is generated by $\cal M.$
\end{theorem}

\bpr
(1) It is easy to see that for each $g_1,...,g_k$ such that $\sum_{i=1}^k v(g_i)=0$,
$\tau_{g_1}\cdot ...\cdot \tau_{g_k}$ is $H^N$ invariant.
Hence, $\tau_{g_1}\cdot ...\cdot \tau_{g_k}\in \C[G^N]^{G\times H^N}$ and, by Proposition \ref{G-G'-free}, it defines a function $\lambda_{g_1,...,g_k}$ in $\C[X_{G/H}(F_N)].$

(2) Since $\tau_g$ generate $\C[G^N]^G,$ for $g\in B,$ every $f\in \C[X_{G/H}(F_N)]=\C[G^N]^{G\times H^N}\subset \C[G^N]^G$
can be written as $$f=\sum_{i=1}^s f_i,$$
where each $f_i$ is a monomial in $\tau_g$, for $g\in B.$
Let $s(f)$ be the smallest such $s$ for given $f.$

We prove the statement of the theorem by contradiction.
Suppose that there is an element of $\C[G^N]^{G\times H^N}$ which is not a polynomial expression in elements of $\cal M.$ Choose such $f$ with $s(f)$ as small as possible.
Let $f=\sum_{i=1}^{s} f_i$ be a corresponding decomposition.
Let $h_1,...,h_N$ be a generating set for $H^N.$ Since each element of $H^N$ acts on each $f_i$ by a scalar multiplication, $h_j\cdot f_i=c_{ji}\cdot f_i$ for some $c_{ji}\in \C^*.$ Then $c_{js}\ne 1$ for some $j$, since otherwise $f_s\in \cal M$ and $f-f_s$ is an element not generated by $\cal M$ with\\ $s(f-f_s)<s=s(f)$.
Since $$f=h_j\cdot f=\sum_{i=1}^{s} c_{ji}f_i,$$
$$c_{js}f-f=c_{js}f-\sum_{i=1}^{s} c_{ji}f_i$$
and, hence,
$$f=\frac{1}{c_{js}-1}\sum_{i=1}^{s(f)-1} (c_{js}-c_{ji})f_i,$$
contradicting the minimality of $s(f).$
\epr

Since $\cal M$ contains $\lambda_{g_1,...,g_k}$ for an arbitrarily long sequences $g_1,...,g_k,$ the set $\cal M$ is usually infinite. We are going to reduce it to a finite generating set now.

Let $\cal V(m,N)$ be the set of all multisets\footnote{A multiset is a ``set" in which an element may appear more than once.} $\{v_1,v_2,...,v_k\}$ of vectors in $(\Z/m)^N,$ such that $\sum_{i=1}^k v_i=0,$ but no proper, non-empty submultiset of $\{v_1,...,v_k\}$ adds up to $0.$
Let $\cal M'\subset \cal M$ be composed of functions $\lambda_{g_1,...,g_k}$ such that $g_1,...,g_k\in B$ and
$\{v(g_1),...,v(g_k)\}\in \cal V(m,N)$.
Since every element of $\cal M$ is a product of elements in $\cal M'$, we have

\begin{corollary} \lb{gen-M'}
$\C[X_{G/H}(F_N)]$ is generated by the elements of $\cal M'.$
\end{corollary}

$\cal V(m,N)$ is finite and, consequently, $\cal M'$ is finite as well.
Indeed, we have:

\begin{proposition} \lb{m^N}
The sequences in $\cal V(m,N)$ have length at most $m^N.$
\end{proposition}

\bpr
Suppose that $\{v_1,v_2,...,v_k\}\in \cal V(m,N).$
Then $\sum_{i=1}^l v_i \ne 0$ in $(\Z/m)^N$, for all $l=1,...,k-1.$ If $k> m^N,$
then $\sum_{i=1}^{l_1} v_i =\sum_{i=1}^{l_2} v_i$ for some $l_1<l_2\leq k$
and $\sum_{i=l_1+1}^{l_2} v_i=0.$
\epr

For example, $\cal V(2,2)$ is composed of five multisets: $\{(0,0)\},$ $\{(1,0),(1,0)\},$ $\{(0,1),(0,1)\},$ $\{(1,1),(1,1)\},$ $\{(1,0),(0,1),(1,1)\}.$

The length of the longest sequence in $\cal V(m,N)$ is called the Davenport constant for the group $(\Z/m)^N.$ No explicit formula for it is known.
However, it is known that
$$N(m-1)+1\leq d(m,N)\leq (m-1)(1+(N-1)m\, ln\, m)+1$$
and, furthermore, the equality on the left holds for $m$ prime (and any $N$) and for $N=2$ (and any $m$), \cite{GLP}.

\bcor 
Since $\C[X_{SL(2,\C)}(F_2)]=\C[\tau_{\gamma_1},\tau_{\gamma_2},\tau_{\gamma_1\gamma_2}]$,
$X_{PSL(2,\C)}(F_2)$ is generated by $g_1=\tau_{\gamma_1}^2,$ $g_2=\tau_{\gamma_2}^2$, $g_3=\tau_{\gamma_1\gamma_2}^2$ and $g_4=\tau_{\gamma_1}\tau_{\gamma_2}\tau_{\gamma_1\gamma_2}.$
It is easy to see that $g_1g_2g_3-g_4^2$ generates all other relations between these generators.
Therefore,
$$X_{PSL(2,\C)}(F_2)=\C[g_1,g_2,g_3,g_4]/(g_1g_2g_3-g_4^2).$$
\ecor

$PSL(2,\C)$-character varieties were studied in further detail in \cite{HP}.

\begin{example}
$SO(4,\C)=(SL(2,\C)\times SL(2,\C))/(\Z/2\Z)$ and $SO(6,\C)=SL(4,\C)/(\Z/2).$ Therefore,
the method of this section provides generating sets for coordinate rings of $SO(4,\C)$- and $SO(6,\C)$-character varieties, alternative to those obtained in \cite[Theorem 8]{S-gen}. We will discuss $SO(4,\C)$-character varieties in further detail in Sec. \ref{s_SO4-desc}.
\end{example}

Corollary \ref{gen-M'} provides finite generating sets for coordinate rings of $G/H$-character varieties for $G=SL(n,\C), Sp(2n,\C),$ and $SO(2n+1,\C).$ (In the last two cases, $H$ is either trivial or $H=\{\pm 1\}$.)

The case of $G=SO(2n,\C)$ is not much different. In this case, the coordinate ring $\C[X_G(\Gamma)]$ is generated by $\tau_g$, for $g$ in some $B\subset \Gamma$, and by $Q_{2n}(g_1,..,g_n),$ for $g_1,...,g_n$ in some $B'\subset \Gamma,$ c.f. \cite[Thm. 8]{S-gen}, where the later are defined as follows:

Let $Q_{2n}: M(2n,\C)^{n}\to \C$ be the ``full polarization" of the function $X\to \text{Pfaffian}(X-X^T)$. It is given explicitly by \vspace*{.1in}

$Q_{2n}(A_1,...,A_n)=\sum_{\sigma\in S_{2n}} sn(\sigma)
 (A_{1,\sigma(1),\sigma(2)}-A_{1,\sigma(2),\sigma(1)})...\hspace*{.2in}$\vspace*{-.2in}\\
$$
\hspace*{2.7in} (A_{n,\sigma(2n-1),\sigma(2n)}-A_{n,\sigma(2n),\sigma(2n-1)}),
$$

where $A_{i,j,k}$ is the $(j,k)$-th entry of the $2n\times 2n$ matrix $A_i$, $i=1,...,n,$ and $sn(\sigma)=\pm 1$ is the sign of $\sigma.$

In \cite{S-gen}, we prove that the function $Hom(\Gamma,SO(2n,\C))\to \C$ sending $\rho$ to $Q_{2n}(\rho(g_1),...,\rho(g_n))$, for some $g_1,...,g_n$, is regular and invariant under the conjugation of $\rho$ by elements of $SO(2n,\C)$. Its factorization to a function on $X_{SO(2n,\C)}(\Gamma)$ is the function $Q_{2n}(g_1,...,g_n)$ alluded to above.

We say that $\tau_g$ has weight $v(g)\in (\Z/2)^N$ and that $Q_{2n}(g_1,..,g_n)$
has weight $v(g_1)+...+v(g_n).$

The only non-trivial central subgroup of $SO(2n,\C)$ is $H=\{\pm 1\}$ and for $PSO(2n,\C)=SO(2n,\C)/H$ we have

\begin{theorem}
$\C[X_{PSO(2n,\C)}(F_N)]$ is generated by monomials in the above generators of $\C[X_{SO(2n,\C)}(F_N)]$ of total weight $0$ in $(\Z/2)^N.$
\end{theorem}

The proof is a straightforward adaptation of that of Theorem \ref{generators-G'}.
\qed


%
\section{An involution $\sigma$ on $X_{SO(2n,\C)}(\Gamma)$}
\lb{s_invol}
%

Let $M$ be any orthogonal hyperplane reflection in $\C^{2n}.$ (Therefore, $M$ is a $2n\times 2n$ orthogonal matrix of determinant $-1$). Then conjugation by $M$ determines a transformation $\sigma$ on $X_{SO(2n,\C)}(\Gamma)$ for any group $\Gamma$ and an induced transformation on $\C[X_{SO(2n,\C)}(\Gamma)]$ which we will denote by the same symbol. Since $M^2\in SO(2n,\C)$, $\sigma$ is an involution.
Furthermore, since every two matrices $M$ (as above) are conjugated one with another by a matrix in $SO(2n,\C)$, the involution $\sigma$ does not depend on the choice of $M.$

\blem \lb{sigmaQ}
$\sigma(Q_{2n}(g_1,...,g_n))=-Q_{2n}(g_1,...,g_n),$ for any $g_1,...,g_n\in \Gamma.$
\elem

\bpr  It suffices to show that for any $2n\times 2n$ matrices $X_1,...,X_n$
$$Q_{2n}(MX_1M^T,...,MX_nM^T)=-Q_{2n}(X_1,...,X_n).$$
Furthermore, by \cite[Prop. 10(2)]{S-gen}, it is enough to assume that $X_1=...=X_n.$
Hence, by \cite[Prop. 10(1)]{S-gen}, this equality reduces to
$$Pf(M(X-X^T)M^T)=-Pf(X-X^T),$$
($Pf$ denotes the Pfaffian), which follows from the fact that
$$Pf(M(X-X^T)M^T)=Pf(X-X^T)\cdot Det(M).$$
\epr

\bcor \lb{QQ-sigma-inv} 
For any $g_1,...,g_n,g_1',...,g_n'\in \Gamma$,
$$Q_{2n}(g_1,...,g_n)Q_{2n}(g_1',...,g_n')\in \C[X_{SO(2n,\C)}(\Gamma)]^\sigma.$$
\ecor


We will denote by $X_{SO(2n,\C)}(\Gamma)/\sigma$ the quotient of $X_{SO(2n,\C)}(\Gamma)$ by the action of $\Z/2=\{1,\sigma\}$ on it. Note that we have an embedding
$$\eta: X_{SO(2n,\C)}(\Gamma)/\sigma \hookrightarrow X_{O(2n,\C)}(\Gamma).$$
(That follows from the fact that any orthogonal matrix is either special orthogonal or a product of special orthogonal with $M$.)

\blem \lb{im-eta}
The image of $\eta$ is isomorphic with $Hom(\Gamma,SO(2n,\C))//O(2n,\C)$.
\elem

\bpr The proof relies on the following statement of invariant theory: For any commutative $\C$-algebra $P$ with an action of a reductive group $G$ on it and for any ideal $I\triangleleft P^G$, the quotient map $\pi: P\to P/IP$ restricts to an isomorphism $\pi^G: P^G/I\to (P/IP)^G$. (Indeed, since $G$ is reductive, $P$ decomposes into a direct sum of $\pi^{-1}(P/IP)$ and of a complementary $G$-module. Therefore, the map $P^G\to (P/IP)^G$ is onto. To see that $\pi^G$ is $1$-$1$, consider any $p\in  P^G$ such that $\pi^G(p+I)=0.$ Then $p=\sum p_j i_j,$ where $p_j\in P$ and $i_j\in I.$ Let $R:P\to P^G$ be the Reynolds operator. Then $p=R(p)=\sum R(p_j i_j)=\sum R(p_j)i_j$ implying that $p\in I.$)

Consider now the function $d_\gamma: X_{SO(2n,\C)}(\Gamma)\to \{\pm 1\}$ sending $[\rho]$ to $det(\rho(\gamma)).$
The image of $\eta$ is the zero set of $I\triangleleft \C[X_{O(2n,\C)}(\Gamma)],$ generated by $d_{\gamma_1},...,d_{\gamma_N},$ for the generators $\gamma_1,...,\gamma_N$ of $\Gamma.$
Now the statement of the lemma follows from the above paragraph for $G=O(2n,\C),$ $P=\C[Hom(\Gamma,G)]$ and the fact that
$P/IP=\C[Hom(\Gamma,SO(2n,\C))].$
\epr

\bpro \lb{eta*}
(1) $\eta: X_{SO(2n,\C)}(\Gamma)/\sigma\to Im\, \eta$ induces an isomorphism\\
$\eta^*: \C[Im\, \eta]\to {\cal T}_{SO(2n,\C)}(\Gamma)\subset\C[X_{SO(2n,\C)}(\Gamma)]^{\sigma}$.\\
(2) For $\Gamma$ free, $\eta: X_{SO(2n,\C)}(\Gamma)/\sigma \hookrightarrow Im\, \eta\subset X_{O(2n,\C)}(\Gamma)$ is an isomorphism.
\epro

\bpr (1) Since $\C[X_{O(2n,\C)}(\Gamma)]$ is generated by trace functions, $\C[Im\, \eta]$ is generated by trace functions as well and, hence, $\eta^*$ maps it to ${\cal T}_{SO(2n,\C)}(\Gamma).$
Since $\eta$ is onto $Im\, \eta$, the dual map, $\eta^*,$ is $1$-$1$. Finally, since every trace function on $X_{SO(2n,\C)}(\Gamma)$ factors through a trace function on $X_{O(2n,\C)}(\Gamma)$, $\eta^*$ is onto.\\
(2) For $\Gamma=\Z$, $\eta$ reduces to an isomorphism between $(SO(2n,\C)//SO(2n,\C))/\sigma$ and $SO(2n,\C)//O(2n,\C)$ (cf. Lemma \ref{im-eta}). For $\Gamma$ of rank $\geq 2$, the result of \cite{Vi} implies that $\eta: X_{SO(2n,\C)}(\Gamma)/\sigma\to Im\, \eta$ is a normalization map. By Lemma \ref{im-eta},
$$Im\, \eta=Hom(\Gamma,SO(2n,\C))//O(2n,\C)$$
which is $SO(2n,\C)^N//O(2n,\C)$ for $\Gamma$ free, and hence it is normal.
\epr

\bcor\lb{cor-T-sigma}
For $\Gamma$ free,
$${\cal T}_{SO(2n,\C)}(\Gamma)=\C[X_{SO(2n,\C)}(\Gamma)]^{\sigma}.$$
\ecor

\bpr follows by comparing images of $\eta^*$ in parts (1) and (2) of Proposition \ref{eta*}.
\epr

%
\section{Trace algebras and full trace algebras for $G=SO(2n,\C)$}
\lb{s_trace_alg}
%

\noindent{\bf Proof of Prop. \ref{module_gen}:}
By \cite[Thm. 8]{S-gen}, $\C[X_{SO(2n,\C)}(\Gamma)]$ is generated by $\tau_\gamma$'s and by functions $Q_{2n}(g_1,...,g_n)$, for words $g_1,...,g_n\in \Gamma$ in $\gamma_1,...,\gamma_N$ of length at most $\nu_n-1$ in which the number of inverses is not larger than the half of the length of the word. Now the statement of Proposition \ref{module_gen} follows from the corollary above and the fact that $\tau_\gamma$'s and products of $Q_{2n}$'s are $\sigma$-invariant (by Corollary \ref{QQ-sigma-inv}).
\qed


\bpro \lb{torsion-over-ft}
For every group $\Gamma,$ the full trace algebra $\FT_{SO(2n,\C)}(\Gamma)$ is a $\T_{SO(2n,\C)}(\Gamma)$-algebra generated by $Q_{2n}(\gamma,...,\gamma)$ for $\gamma\in \Gamma.$
\epro

\bpr
By \cite[Sec 23.2]{FH}, the representation ring of $SO(2n,\C)$ is generated by the exterior powers, $\wedge^k V,$ of the defining ($2n$-dimensional) representation $V$ and by the representations $D_n^+,D_n^-,$ whose highest weights are twice that of the $\pm$-half-spin representations. (The last two representations are denoted by $D_{\pm}$ in \cite{S-gen}.)

For every $M\in M(2n,\C)$ with eigenvalues $\lambda_1,...,\lambda_{2n},$ the induced transformation on $\wedge^k \C^{2n}$ has trace $\sum_{i_1<...<i_k} \lambda_{i_1}\cdot ...\cdot \lambda_{i_k}$, which is a polynomial in $tr(M^d)$ for $d=1,2,3,...$ That polynomial depends on $k$ only.
Therefore, for every $k$, $\tau_{\gamma,\wedge^k V}\in \T_{SO(2n,\C)}(\Gamma)$, implying that  $\FT_{SO(2n,\C)}(\Gamma)$ is a $\T_{SO(2n,\C)}(\Gamma)$-algebra generated by $\tau_{\gamma,D_n^\pm},$ for $\gamma\in \Gamma.$ Since
$D_n^++D_n^-$ is a polynomial in $V$, cf. \cite[Sec 23.2]{FH},
$$\tau_{\gamma,D_n^+}+\tau_{\gamma,D_n^-}\in \T_{SO(2n,\C)}(\Gamma)$$
and, consequently, $\FT_{SO(2n,\C)}(\Gamma)$ is generated by expressions $\tau_{\gamma,D_n^+}-\tau_{\gamma,D_n^-}$ over $\T_{SO(2n,\C)}(\Gamma)$. By \cite[Prop 10]{S-gen},
$$\tau_{\gamma,D_n^+}-\tau_{\gamma,D_n^+}=(2i)^{-n}(n!)^{-1}Q_{2n}(\gamma,...,\gamma),$$
implying
that $Q_{2n}(\gamma,...,\gamma)$, for $\gamma\in \Gamma,$ are $\T_{SO(2n,\C)}(\Gamma)$-algebra generators of $\FT_{SO(2n,\C)}(\Gamma).$
\epr

We analyze the extensions
$$\T_{SO(2n,\C)}(\Gamma)\subset \FT_{SO(2n,\C)}(\Gamma)\subset \C[X_{SO(2n,\C)}(\Gamma)]$$
in further detail for $\Gamma=F_2$ and $n=2$ in Sections \ref{s_SO4-desc}-\ref{s_Coord-gen2}. In particular, we prove Theorem \ref{Coord-gen} there.

%
\section{A presentation of $\C[X_{SO(4,\C)}(F_2)]$ in terms of generators and relations}
\label{s_SO4-desc}
%

There is a natural isomorphism $SL(2,\C)\times SL(2,\C)\to Spin(4,\C)$ inducing the epimorphism
\beq\lb{e_isom}
\phi: SL(2,\C)\times SL(2,\C) \to SO(4,\C)
\eeq
with the kernel $\{(I,I),(-I,-I)\},$ \cite[Ch. 3.\S 2 Example 2]{GOV}.
This epimorphism can be defined as follows: Consider the action of $SL(2,\C)\times SL(2,\C)$ on
$\C^2\otimes \C^2=\C^4$ by matrix multiplication.
That action preserves the symmetric, non-degenerate product
$$((u_1,u_2),(v_1,v_2))=det(u_1,v_1)det(u_2,v_2),$$
where $det(u,v)$ is the determinant of the $2\times 2$ matrix composed of vectors $u$ and $v.$
Consequently, it induces the desired epimorphism $\phi.$

Consider the following orthonormal basis of $\C^2\otimes \C^2:$
$$w_1=\frac{1}{\sqrt{2}}(e_1\otimes e_1+e_2\otimes e_2),\quad
w_2=\frac{i}{\sqrt{2}}(e_1\otimes e_1-e_2\otimes e_2),$$
$$w_3=\frac{i}{\sqrt{2}}(e_1\otimes e_2+e_2\otimes e_1),\quad
w_4=\frac{1}{\sqrt{2}}(e_1\otimes e_2-e_2\otimes e_1),$$
where $e_1=(1,0),$ $e_2=(0,1).$
A direct computation (which will be useful later) shows that for
$$((a_{ij}), (b_{ij}))\in SL(2,\C)\times SL(2,\C),$$
$\phi((a_{ij}), (b_{ij}))$ is given by the following matrix with respect to the above basis:
\beq\lb{e_phi}
\phi((a_{ij}), (b_{ij}))=
\eeq
$\frac{1}{2}\cdot\left(
{\scriptsize
\begin{array}{cccc}
a_{11}b_{11}+a_{12}b_{12}+a_{21}b_{21}+a_{22}b_{22} &
ia_{11}b_{11}-ia_{12}b_{12}+ia_{21}b_{21}-ia_{22}b_{22} \\
-ia_{11}b_{11}-ia_{12}b_{12}+ia_{21}b_{21}+ia_{22}b_{22} &
a_{11}b_{11}-a_{12}b_{12}-a_{21}b_{21}+a_{22}b_{22} \\
-ia_{21}b_{11}-ia_{22}b_{12}-ia_{11}b_{21}-ia_{12}b_{22} &
a_{21}b_{11}-a_{22}b_{12}+a_{11}b_{21}-a_{12}b_{22} \\
-a_{21}b_{11}-a_{22}b_{12}+a_{11}b_{21}+a_{12}b_{22} &
-ia_{21}b_{11}+ia_{22}b_{12}+ia_{11}b_{21}-ia_{12}b_{22}
\end{array}
}
\right.$\\
\hspace*{.7in}$\left.
{\scriptsize
\begin{array}{cccc}
ia_{12}b_{11}+ia_{11}b_{12}+ia_{22}b_{21}+ia_{21}b_{22} &
-a_{12}b_{11}+a_{11}b_{12}-a_{22}b_{21}+a_{21}b_{22} \\
a_{12}b_{11}+a_{11}b_{12}-a_{22}b_{21}-a_{21}b_{22} &
ia_{12}b_{11}-ia_{11}b_{12}-ia_{22}b_{21}+ia_{21}b_{22} \\
a_{22}b_{11}+a_{21}b_{12}+a_{12}b_{21}+a_{11}b_{22} &
ia_{22}b_{11}-ia_{21}b_{12}+ia_{12}b_{21}-ia_{11}b_{22} \\
-ia_{22}b_{11}-ia_{21}b_{12}+ia_{12}b_{21}+ia_{11}b_{22} &
a_{22}b_{11}-a_{21}b_{12}-a_{12}b_{21}+a_{11}b_{22}
\end{array}
}
\right).$\vspace*{.15in}

We are going to use the above description of $SO(4,\C)$ and the method of Sec. \ref{s_G-G'} to find a presentation of $\C[X_{SO(4,\C)}(F_2)]$ in terms of generators and relations.

By Proposition \ref{G-G'-free}, $$X_{SO(4,\C)}(F_2)=(X_{SL(2,\C)}(F_2)\times X_{SL(2,\C)}(F_2))/(\Z/2\times \Z/2).$$
The above action of $\Z/2\times \Z/2$ on $X_{SL(2,\C)}(F_2)\times X_{SL(2,\C)}(F_2)$ can be described as follows: $(\ve_1,\ve_2)\in \{\pm 1\}\times \{\pm 1\}=\Z/2\times \Z/2$ sends the equivalence class of
$\rho:\Gamma\to SL(2,\C)\times SL(2,\C)$ to the equivalence class of $\rho'$ such that $$\rho'(\gamma_i)=(\ve_i\rho_1(\gamma_i),\ve_i\rho_2(\gamma_i))$$ for $i=1,2.$


Let us abbreviate the generators $\tau_{\gamma_1}, \tau_{\gamma_2}, \tau_{\gamma_1\gamma_2}$ of $\C[X_{SL(2,\C)}(F_2)]$ by $\tau_1,\tau_2,\tau_{12}.$ (Hence, $\C[X_{SL(2,\C)}(F_2)]=\C[\tau_1,\tau_2,\tau_{12}].$) We will denote the generators of
$$\C[X_{SL(2,\C)}(F_2)\times X_{SL(2,\C)}(F_2)]=\C[X_{SL(2,\C)}(F_2)]\otimes \C[X_{SL(2,\C)}(F_2)]$$
by $\tau_{1,i},\tau_{2,i},\tau_{12,i},$ where $i=1,2$ indicates the first or the second copy of $\C[X_{SL(2,\C)}(F_2)]$.

\bcor \lb{SO4-gen}
$$\C[X_{SO(4,\C)}(F_2)]=\C[X_{SL(2,\C)}(F_2)\times X_{SL(2,\C)}(F_2)]^{\Z/2\times \Z/2}$$
is generated by
$$a_{i,j,k}=\tau_{i,j}\tau_{i,k},\quad b_{j,k}=\tau_{12,j}\tau_{12,k},\quad \text{for $i,j,k=1,2$ where $j\leq k$,}$$
and by
$$c_{i,j,k}=\tau_{1,i}\tau_{2,j}\tau_{12,k}, \quad \text{for $i,j,k=1,2$}.$$
\ecor

\bpro \lb{SO4-desc}
$\C[X_{SO(4,\C)}(F_2)]$ is the quotient of the polynomial ring in the above $17$ generators $a_{i,j,k},b_{i,j},c_{i,j,k}$ by the ideal generated by
$$a_{i,j,k}a_{i,j',k'}-a_{i,j,j'}a_{i,k,k'},\quad b_{j,k}b_{j',k'}-b_{j,j'}b_{k,k'},\quad c_{i,j,k}\cdot c_{i',j',k'}-a_{1,i,i'}a_{2,j,j'}b_{k,k'},$$
$$a_{1,j,k}c_{i,j',k'}-a_{1,i,k}c_{j,j',k'},\quad a_{2,j,k}c_{i,j',k'}-a_{2,j',k}c_{i,j,k'},\quad
b_{j,k}c_{i,j',k'}-b_{j,k'}c_{i,j',k},$$
where $a_{i,2,1}=a_{i,1,2}$ and $b_{2,1}=b_{1,2}.$
\epro

\bpr
Since all $\tau_{i,j}$'s and $\tau_{12,i}$'s are algebraically independent, all relations between the generators of $\C[X_{SO(4,\C)}(F_2)]$ follow from different decompositions of monomials in $\tau$'s into products of $a$'s, $b$'s, and $c$'s. It is straightforward to check that any two such decompositions are related by the relations listed above.
\epr
%
\section{Proof of Theorem \ref{Coord-gen}(1) and (2)}
\lb{s_Coord-gen1}
%

Consider the isomorphism $\phi: SL(2,\C)\times SL(2,\C))/\{\pm (I,I)\} \to SO(4,\C)$ of Sec. \ref{s_SO4-desc}.

\blem \lb{sigma-phi}
If $M=\left({\scriptsize \be{matrix} 1&0&0&0\\ 0&1&0&0\\ 0&0&1&0\\ 0&0&0&-1 \en{matrix} }\right),$
then the involution $\sigma$ on $SO(4,\C),$ $\sigma(X)=M\cdot X\cdot M^{-1},$ satisfies $\sigma \phi(A,B)=\phi(B,A)$ for any $A,B\in SL(2,\C).$
\elem

\bpr
Since
$$C_1=\left(\be{matrix} 1 & 1\\ 0 & 1\en{matrix}\right)\ \text{and}\ C_2=\left(\be{matrix} 0 & -1\\ 1 & 0\en{matrix}\right)$$ generate $SL(2,\Z)$, which is a Zariski-dense subgroup of $SL(2,\C)$,
and since $\sigma$ is an algebraic group automorphism, it is enough to verify the statement
for $(A,B)=(C_i,C_j),$ $i,j=1,2.$
A substitution of $(C_i,C_j)$ into (\ref{e_phi}) yields\\
$$\phi(C_1,C_1)=\frac{1}{2}\left(
{\scriptsize
\begin{array}{cccc}
 3 & -i & 2i & 0 \\
 -i & 1 & 2 & 0 \\
 -2i & -2 & 2 & 0 \\
 0 & 0 & 0 & 2 \\
\end{array}
}
\right),\
\phi(C_1,C_2)=\frac{1}{2}\left(
{\scriptsize
\begin{array}{cccc}
-1 & i & 0 & -2 \\
i & 1 & -2 & 0 \\
0 & 2 & 1 & i \\
2 & 0 & i & -1 \\
\end{array}
}
\right),$$
$$\phi(C_2,C_1)=\frac{1}{2}\left(
{\scriptsize
\begin{array}{cccc}
 -1 & i & 0 & 2 \\
 i & 1 & -2 & 0 \\
 0 & 2 & 1 & -i \\
 -2 & 0 & -i & -1 \\
\end{array}
}
\right),\
\phi(C_2,C_2)=\left(
{\scriptsize
\be{matrix} 1&0&0&0\\ 0&-1&0&0\\ 0&0&-1&0\\ 0&0&0&1 \en{matrix}
}
\right).$$
Since
$$\sigma\left(
{\scriptsize
\begin{array}{cccc}
 x_{11} & x_{12} & x_{13} & x_{14} \\
 x_{21} & x_{22} & x_{23} & x_{24} \\
 x_{31} & x_{32} & x_{33} & x_{34} \\
 x_{41} & x_{42} & x_{43} & x_{44} \\
\end{array}
}
\right)\ =\left(
{\scriptsize
\begin{array}{cccc}
x_{11} & x_{12} & x_{13} & -x_{14} \\
 x_{21} & x_{22} & x_{23} & -x_{24} \\
 x_{31} & x_{32} & x_{33} & -x_{34} \\
 -x_{41} & -x_{42} & -x_{43} & x_{44} \\
\end{array}
}
\right),$$
we easily see that $\phi(C_1,C_1)$ and $\phi(C_2,C_2)$ are $\sigma$-invariant,
while\\
$\sigma(\phi(C_1,C_2))=\phi(C_2,C_1).$
\epr

\bcor \lb{sigma-on-gen}
$\sigma(a_{i,j,k})=a_{i,\hat j,\hat k},$ $\sigma(b_{j,k})=b_{\hat j,\hat k},$
$\sigma(c_{i,j,k})=c_{\hat i,\hat j,\hat k},$
where $\hat 1=2$ and $\hat 2=1.$
\ecor

\blem
If $\sigma$ is an involution on a commutative $\C$-algebra $R$ generated by
$r_1,...,r_s$, then\\
(1) $R^\sigma$ is generated by $p_i=r_i+\sigma(r_i)$ and by $q_iq_j,$ for $i,j=1,...,s,$ where
$q_i=r_i-\sigma(r_i).$\\
(2) As an $R^\sigma$-module, $R$ is generated by $1$ and by $q_i$ for $i=1,...,s.$
\elem

\bpr (1) Since $p_i$'s and $q_i$'s form a generating set of $R,$ every element $r\in R^\sigma$ is of the form $r=\sum_k m_k$, where $m_j$'s are monomials in $p_i$'s and $q_i$'s. Clearly $\sigma m=\pm m$ for every monomial $m$ in that polynomial.
Since $r=\frac{r+\sigma r}{2},$ by replacing every $m_k$ by $\frac{m_k+\sigma m_k}{2}$ if necessary, we can assume without loss of generality that all monomials, $m_k$, are $\sigma$-invariant. Therefore, they are products of $p_i$'s and $q_iq_j$'s.

(2) Clearly $R$ is generated by $q_i$'s as an $R^\sigma$-algebra. Since $q_iq_j\in R^\sigma,$
the algebra $R$ is generated by $1$ and by $q_i$'s as an $R^\sigma$-module.
\epr

By Corollary \ref{sigma-on-gen} and by the above lemma, we have

\bcor\lb{cor-rs}
(1) As an $\C[X_{SO(4,\C)}(F_2)]^\sigma$-module, $\C[X_{SO(4,\C)}(F_2)]$ is generated by
$$t_0=1,\ t_1=c_{1,1,2}-c_{2,2,1},\ t_2=c_{1,2,1}-c_{2,1,2},\ t_3=c_{1,2,2}-c_{2,1,1},$$
$$t_4=a_{1,1,1}-a_{1,2,2},\ t_5=a_{2,1,1}-a_{2,2,2},\ t_6=b_{1,1}-b_{2,2},\ t_7=c_{1,1,1}-c_{2,2,2}.$$
(This specific order of $t_i$'s will prove convenient later.)\\
(2) As a $\C$-algebra, $\C[X_{SO(4,\C)}(F_2)]^\sigma$ is generated by
$$a_{i,1,1}+a_{i,2,2},\  a_{i,1,2},\ \text{for\ } i=1,2,\  b_{1,2},\  b_{1,1}+b_{2,2},$$
$$c_{1,1,1}+c_{2,2,2},\ c_{1,1,2}+c_{2,2,1},\ c_{1,2,1}+c_{2,1,2},\ c_{1,2,2}+c_{2,1,1},$$
and by the products $t_jt_k,$ for $j,k=1,...,7$
\ecor

\blem 
$Q_4(g_1,g_2)=4(\tau_{g_1,2}\tau_{g_2,2}\tau_{g_1g_2,1}-\tau_{g_1,1}\tau_{g_2,1}\tau_{g_1g_2,2})$\\
in $\C[X_{SO(4,\C)}(F_2)].$
\elem

\bpr A computer algebra system computation based on the definition of $Q_4(\cdot,\cdot)$ and on formula (\ref{e_phi}).
\epr

Now we can conclude the proofs of Theorem \ref{Coord-gen}(1) and (2).

In what follows, we use two classical identities:
$$\tau_{g^{-1}}=\tau_{g}\ \text{and}\ \tau_{g}\tau_{h}=\tau_{gh}+\tau_{gh^{-1}}$$
for $\tau_g,\tau_h \in \C[X_{SL(2,\C)}(\Gamma)].$

By the lemma above,
$$Q_4(g,g)=4(\tau_{g,2}^2\tau_{g^2,1}-\tau_{g,1}^2\tau_{g^2,2})=
4(\tau_{g,2}^2(\tau_{g,1}^2-2)-\tau_{g,1}^2(\tau_{g,2}^2-2))=8(\tau_{g,1}^2-\tau_{g,2}^2).$$
Consequently,
\beq\lb{e_t123}
t_1=-\frac{1}{4}Q_4(\gamma_1,\gamma_2),\ t_4=\frac{1}{8}Q_4(\gamma_1,\gamma_1),\ t_5=\frac{1}{8}Q_4(\gamma_2,\gamma_2),\ t_6=\frac{1}{8}Q_4(\gamma_1\gamma_2,\gamma_1\gamma_2).
\eeq
Since
$$\frac{1}{4}Q_4(\gamma_2\gamma_1^{-1},\gamma_1)=
\tau_{\gamma_2\gamma_1^{-1},2}\tau_{1,2}\tau_{2,1}-\tau_{\gamma_2\gamma_1^{-1},1}\tau_{1,1}\tau_{2,2}=$$
$$(\tau_{1,2}\tau_{2,2}-\tau_{12,2})\tau_{1,2}\tau_{2,1}-(\tau_{1,1}\tau_{2,1}-\tau_{12,1})\tau_{1,1}\tau_{2,2}=$$
$$\tau_{2,1}\tau_{2,2}(\tau_{1,2}^2-\tau_{1,1}^2)-(\tau_{1,2}\tau_{2,1}\tau_{12,2}-\tau_{1,1}\tau_{2,2}\tau_{12,1})=
-a_{2,1,2}t_4+t_2,$$
we have
\beq\lb{e_t2}
t_2=\frac{1}{4}Q_4(\gamma_2\gamma_1^{-1},\gamma_1)+a_{2,1,2}t_4.
\eeq
Analogously,
\beq\lb{e_t3}
t_3=-\frac{1}{4}Q_4(\gamma_1\gamma_2^{-1},\gamma_2)-a_{1,1,2}t_5
\eeq
and, finally, similar computations show that
$$t_7=\frac{1}{4}\left(t_4(a_{2,1,1}+a_{2,2,2})+t_5(a_{1,1,1}+a_{1,2,2})\right)+\frac{1}{2}t_6-
\frac{1}{16}Q_4(\gamma_1\gamma_2^{-1},\gamma_1\gamma_2^{-1}).$$


Now, Corollaries \ref{cor-T-sigma} and \ref{cor-rs}(1) together with these identities implies part (1) of Theorem \ref{Coord-gen}.
Part (2) follows immediately from part (1) and from Proposition \ref{torsion-over-ft}.\qed

%
\section{Proof of Theorem \ref{Coord-gen}(3)}
\lb{s_Coord-gen2}
%

We proceed the proof with several preliminaries.

\comment{By \cite[\S 23.2]{FH}, the representation ring of $SO(4,\C)$ is generated by the $4$-dimensional defining representation, which we will denote by $\C^4$ and by two $3$-dimensional representations, which following \cite[\S 23.2]{FH}, we denote by $D_2^\pm$. (The highest weights of $D_2^\pm$
are twice that of the half-spin representations of $so(4,\C).$)\marginpar{Remove these 2 par}

Since
$$Tr(\rho_1(\gamma)\oplus \rho_2(\gamma))=Tr(\rho_1(\gamma))+ Tr(\rho_2(\gamma)),$$
$$Tr(\rho_1(\gamma)\otimes \rho_2(\gamma))=Tr(\rho_1(\gamma))\cdot Tr(\rho_2(\gamma))$$
for any representations $\rho_1, \rho_2$ of $SO(4,\C)$
and since every representation of $SO(4,\C)$ is a sum of products of the defining representations and of $D_2^\pm$'s, the full trace algebra,
$\FT_{SO(4,\C)}(F_2),$ is generated by $\tau_{\gamma,\C^4}$ and $\tau_{\gamma,D_2^\pm}$ for $\gamma\in F_2.$
}

\blem \lb{characters}
The following identities hold in $\C[X_{SO(4,\C)}(F_2)]$ for every\\ 
$\gamma\in F_2$:\\
(1) $\tau_{\gamma}=\tau_{\gamma,1}\cdot\tau_{\gamma,2}$. (As before, $\tau_{\gamma}$ is an abbreviation for $\tau_{\gamma,\C^4},$ where $\C^4$ is the defining representation of $SO(4,\C).$)\\
(2) One of $\tau_{\gamma,D_2^+}, \tau_{\gamma,D_2^-}$ equals to $\tau_{\gamma,1}^2-1$ and the other to $\tau_{\gamma,2}^2-1$.
\elem

\bpr
(1) Let $\pi_i: SL(2,\C)\times SL(2,\C)\to SL(2,\C)$ by the projection on the $i$-th factor, for $i=1,2.$ Then $\pi_1\otimes \pi_2$ is an irreducible $4$-dimensional representation of $SL(2,\C)\times SL(2,\C)$. Since it sends $(-I,-I)$ to the identity, it factors to
an irreducible $4$-dimensional representation of $SO(4,\C).$ Since the defining representation is the only irreducible $4$-dimensional representation of $SO(4,\C)$ (up to conjugation),
$Tr\pi_1\otimes\pi_2(A_1,A_2)=Tr(\phi(A_1,A_2))$ for every $A_1,A_2\in SL(2,\C)$, where $\phi$ is the representation (\ref{e_isom}). Since $$Tr\pi_1\otimes\pi_2(A_1,A_2)=Tr A_1\cdot Tr A_2,$$ the statement follows.

(2) Let $\eta$ be an irreducible $3$-dimensional representation of $SL(2,\C).$ ($\eta$ is unique up to conjugation.) It is easy to see that $\eta\pi_1$ and $\eta\pi_2$ factor to two nonequivalent irreducible representations of
$$SO(4,\C)=(SL(2,\C)\times SL(2,\C))/{\pm (I,I)}.$$
By the classification of representations of $SO(4,\C)$ (as described in the proof of Proposition \ref{torsion-over-ft}) the only irreducible $3$-dimensional representations of $SO(4,\C)$ are $D_2^+$ and $D_2^-.$ Since $Tr \eta(A)=(Tr A)^2-1,$ the statement follows.
\epr

\bpro \lb{FT-gen}
(1) $\tau_{i,j}^2, \tau_{12,i}^2, \tau_{i,1}\tau_{i,2}, \tau_{12,1}\tau_{12,2},$ $\tau_{1,i}\tau_{2,i}\tau_{12,i},$ for $i,j\in \{1,2\},$
$$\tau_{1,1}\tau_{2,2}\tau_{12,1}+\tau_{1,2}\tau_{2,1}\tau_{12,2},\quad
\tau_{1,1}\tau_{2,2}\tau_{12,2}+\tau_{1,2}\tau_{2,1}\tau_{12,1},$$
and
$$\tau_{1,1}\tau_{2,1}\tau_{12,2}+\tau_{1,2}\tau_{2,2}\tau_{12,1}$$
belong to $\FT_{SO(4,\C)}(F_2)$.\\
(2) $\FT_{SO(4,\C)}(F_2)$ is generated by the above elements as a $\C$-algebra.
\epro

\noindent {Proof of Proposition \ref{FT-gen}(1):} $\tau_{i,j}^2, \tau_{12,i}^2,\tau_{i,1}\tau_{i,2}, \tau_{12,1}\tau_{12,2} \in \FT_{SO(4,\C)}(F_2)$ by Lemma \ref{characters}.

By squaring both sides of $$\tau_{\gamma_1^2\gamma_2}=\tau_1\tau_{12}-\tau_2$$
we obtain
$$\tau_1\tau_2\tau_{12}=-\frac{1}{2}\left(\tau_{\gamma_1^2\gamma_2}^2-\tau_1^2\tau_{12}^2- \tau_2^2\right)$$
in $\C[X_{SL(2,\C)}(F_2)].$
Therefore, $\tau_{1,i}\tau_{2,i}\tau_{12,i}\in \FT_{SO(4,\C)}(F_2),$ for $i=1,2,$ by the Lemma \ref{characters}(2).
Similarly, by applying $\gamma=\gamma_1^2\gamma_2$ to Lemma \ref{characters}(1), we see that
$$\tau_{1,1}\tau_{2,2}\tau_{12,1}+\tau_{1,2}\tau_{2,1}\tau_{12,2}=
\tau_{1,1}\tau_{1,2}\tau_{12,1}\tau_{12,2}+\tau_{2,1}\tau_{2,2}-
\tau_{\gamma_1^2\gamma_2,1}\tau_{\gamma_1^2\gamma_2,2}$$
belongs to $\FT_{SO(4,\C)}(F_2)$ as well.
In the same vein, we prove that
$$\tau_{1,1}\tau_{2,2}\tau_{12,2}+\tau_{1,2}\tau_{2,1}\tau_{12,1},\  \tau_{1,1}\tau_{2,1}\tau_{12,2}+\tau_{1,2}\tau_{2,2}\tau_{12,1}\in \FT_{SO(4,\C)}(F_2)$$
by taking $\gamma=\gamma_1\gamma_2^2$ and $\gamma=\gamma_1^{-1}\gamma_2,$ respectively.\qed\vspace*{.15in}

In order to prove part (2) of Proposition \ref{FT-gen}
consider a $\Z/2\times \Z/2$ action on $X_{SL(2,\C)}(F_2)$ defined for $(k_1,k_2)\in \Z/2\times \Z/2$
by $(k_1,k_2)\cdot [\rho]=[\rho'],$
where $\rho'(\gamma_i)=(-1)^{k_i}\rho(\gamma_i),$ for $i=1,2.$ That action defines a $\Z/2\times \Z/2$-grading on $R=\C[X_{SL(2,\C)}(F_2)]$, whose homogeneous components we denote by $R_{i,j},$ $i,j\in \Z/2.$ More concretely, the elements $f\in R_{i,j},$ for $i,j\in \{0,1\},$ have the property that $(1,0)\cdot f=(-1)^i f$ and $(0,1)\cdot f=(-1)^j f$ for $i,j\in \Z/2.$
In particular, $deg(\tau_1)=(1,0)$, $deg(\tau_2)=(0,1)$ and $deg(\tau_{12})=(1,1).$

The following is straightforward:

\blem\lb{lem-gens}
$R_{0,0}=\C[\tau_1^2,\tau_2^2,\tau_{12}^2,\tau_1\tau_2\tau_{12}]$ and, as $R_{0,0}$-modules,\\
$R_{1,0}$ is generated by $\tau_1$ and $\tau_2\tau_{12},$\\
$R_{0,1}$ is generated by $\tau_2$ and $\tau_1\tau_{12}$, and\\
$R_{1,1}$ is generated by $\tau_{12}$ and $\tau_1\tau_2$.
\elem

Let $R_{0,0,i}$ be the subalgebra of $\C[X_{SL(2,\C)}(F_2)\times X_{SL(2,\C)}(F_2)]$ generated by the generators of $R_{0,0}$ with the second index $i$: $\tau_{1,i}^2,\tau_{2,i}^2, \tau_{12,i}^2, \tau_{1,i}\tau_{2,i}\tau_{12,i}$ for $i=1,2.$
Similarly, let $R_{1,0,k}$ be an $R_{0,0,k}$-module generated by $\tau_{1,k}$ and $\tau_{2,k}\tau_{12,k},$
let $R_{0,1,k}$ be an $R_{0,0,k}$-module generated by $\tau_{2,k}$ and $\tau_{1,k}\tau_{12,k},$ and
let $R_{1,1,k}$ be an $R_{0,0,k}$-module generated by $\tau_{12,k}$ and $\tau_{1,k}\tau_{2,k}.$
\vspace*{.2in}

\noindent {Proof of Proposition \ref{FT-gen}(2):} We need to prove that $\FT_{SO(4,\C)}(F_2) \subset \cal S$, where $\cal S$ denotes the $\C$-subalgebra of $\C[X_{SL(2,\C)}(F_2)\times X_{SL(2,\C)}(F_2)]$ generated by the elements listed in part (1).

Proof that $\tau_{\gamma,D_2^\pm}\in \cal S$ for every $\gamma\in F_2:$ By Lemma \ref{characters}(1), $\tau_{\gamma,D_2^\pm}=\tau_{\gamma,i}^2-1$ for $i=1$ or $2.$ Since $\tau_{\gamma}^2$ is invariant under the $\Z/2\times \Z/2$ action, $\tau_\gamma^2\in R_{0,0}$ and $\tau_{\gamma,D_2^\pm}\in R_{0,0,i}$.
Since $R_{0,0,1}$ and $R_{0,0,2}$ are subalgebras of $\cal S,$ the statement follows.

Proof that $\tau_{\gamma}\in \cal S:$
Consider the epimorphism $\eta:F_2\to \Z/2\times \Z/2$ sending $\gamma_1$ to $(1,0)$ and $\gamma_2$ to $(0,1).$ Since $(1,0)\cdot \tau_\gamma\in \C[X_{SL(2,\C)}(F_2)]$ is either $\tau_\gamma$ or $-\tau_\gamma$ depending on the first component of $\eta(\gamma)$ and $(0,1)\cdot \tau_\gamma$ is either $\tau_\gamma$ or $-\tau_\gamma$ depending on the second component of $\eta(\gamma)$, we see that $\tau_\gamma\in R_{\eta(\gamma)}$ for every $\gamma\in F_2.$

For every $\gamma\in F_2$, consider a three-variable polynomial $p_\gamma$ such that $\tau_\gamma=p(\tau_1,\tau_2,\tau_{12})$.
By Lemma \ref{characters}(1),
$\tau_{\gamma,\C^4}=\tau_{\gamma,1}\cdot \tau_{\gamma,2}.$
The following lemma completes the statement of Proposition \ref{FT-gen}(2).
\qed

\blem 
For every $p$ as above,
$p(\tau_{1,1},\tau_{2,1},\tau_{12,1})\cdot p(\tau_{1,2},\tau_{2,2},\tau_{12,2})\in \cal S.$
\elem

\bpr We consider the four possible values of $\eta(\gamma).$
If $\eta(\gamma)=(0,0)$, then the statement is obvious since $R_{0,0,1},R_{0,0,2}\subset \cal S.$
If $\eta(\gamma)=(1,0)$, then $p(\tau_{1,1},\tau_{2,1}\tau_{12,1})\cdot p(\tau_{1,2},\tau_{2,2}\tau_{12,2})$ belongs to an $S$-module generated by
$$\tau_{1,1}\tau_{1,2},\ \tau_{2,1}\tau_{12,1}\tau_{2,2}\tau_{12,2}\ \text{and } \tau_{1,1}\tau_{2,2}\tau_{12,2}+\tau_{1,2}\tau_{2,1}\tau_{12,1}.$$
All these elements belong to $\cal S.$
One can perform a similar verification for $\eta(\gamma)=(0,1)$ and for $(1,1).$
\epr

\noindent {\bf Proof of Theorem \ref{Coord-gen}(3):}
Define a $\Z_{\geq 0}$-grading on $\C[\tau_{i,j},\tau_{12,j},\ i,j=1,2]$ by declaring that all $\tau_{i,j}$'s and $\tau_{12,j}$'s have degree $1.$ Note that $\C[X_{SO(4,\C)}(F_2)]$ is a graded subalgebra of $\C[\tau_{i,j},\tau_{12,j},\ i,j=1,2]$. Since the generators of $\FT_{SO(4,\C)}(F_2)$ of Proposition \ref{FT-gen} are homogeneous, $\FT_{SO(4,\C)}(F_2)$ is in turn a graded subalgebra of $\C[X_{SO(4,\C)}(F_2)].$

Suppose $Q_4(\gamma_2\gamma_1^{-1},\gamma_1)$ is a $\FT_{SO(4,\C)}(F_2)$-linear combination of $1,$ $Q_4(\gamma_1,\gamma_2),$ $Q_4(\gamma_1\gamma_2^{-1},\gamma_2).$ Then from equalities (\ref{e_t123}), (\ref{e_t2}), (\ref{e_t3}), we see that $t_3$ is a\\
$\FT_{SO(4,\C)}(F_2)$-linear combination of $t_0,t_1,t_2.$
Since 
$$deg(t_0)=0,\quad deg(t_1)=deg(t_2)=deg(t_3)=3,$$ 
and $\FT_{SO(4,\C)}(F_2)$ is graded,
$$t_3=c_0t_0+c_1t_1+c_2t_2,$$
for some $c_0,c_1,c_2\in \FT_{SO(4,\C)}(F_2)$ such that 
$$deg(c_0)=3\ \text{and}\ deg(c_1)=deg(c_2)=0.$$ 
Therefore, a non-trivial $\C$-linear combination of $t_1,t_2,t_3$ belongs to
$\FT_{SO(4,\C)}(F_2)$. Since such linear combination is a homogeneous element of degree $3$, it has to coincide with a $\C$-linear combination of the generators of $\C[X_{SO(4,\C)}(F_2)]$ of degree $3$ listed in Proposition \ref{FT-gen}. It is easy to see that it is impossible, since $\tau_{i,j}$'s and $\tau_{12,j}$'s are algebraically independent.

The proof of the non-redundancy of the other three generators of Theorem \ref{Coord-gen}(2) is analogous.
\qed

%

\end{document}